\documentclass[12pt,reqno]{article}

\usepackage[usenames]{color}
\usepackage{amssymb}
\usepackage{graphicx}
\usepackage{amscd}

\usepackage[colorlinks=true,
linkcolor=webgreen,
filecolor=webbrown,
citecolor=webgreen]{hyperref}

\definecolor{webgreen}{rgb}{0,.5,0}
\definecolor{webbrown}{rgb}{.6,0,0}

\usepackage{algorithm}
\usepackage{algpseudocode}
\algnewcommand{\algorithmicand}{\textbf{ and }}
\algnewcommand{\algorithmicor}{\textbf{ or }}
\algnewcommand{\OR}{\algorithmicor}
\algnewcommand{\AND}{\algorithmicand}
\algnewcommand{\var}{\texttt}

\usepackage{color}
\usepackage{fullpage}
\usepackage{float}

\usepackage{graphics,amsmath,amssymb}
\usepackage{amsthm}
\usepackage{amsfonts}
\usepackage{latexsym}
\usepackage{epsf}

\newcommand{\seqnum}[1]{\href{http://oeis.org/#1}{\underline{#1}}}
\begin{document}
	
	\theoremstyle{plain}
	\newtheorem{theorem}{Theorem}
	\newtheorem{corollary}[theorem]{Corollary}
	\newtheorem{lemma}[theorem]{Lemma}
	\newtheorem{proposition}[theorem]{Proposition}
	
	\theoremstyle{definition}
	\newtheorem{definition}[theorem]{Definition}
	\newtheorem{example}[theorem]{Example}
	\newtheorem{conjecture}[theorem]{Conjecture}
	
	\theoremstyle{remark}
	\newtheorem{remark}[theorem]{Remark}	
	
	\begin{center}
		\vskip 1cm{\LARGE\bf On the Density of Spoof Odd Perfect Numbers}
		\vskip 1cm
		\large
		L\'aszl\'o T\'oth\\
		Rue des Tanneurs 7 \\
		L-6790 Grevenmacher \\
		Grand Duchy of Luxembourg \\
		\href{mailto:uk.laszlo.toth@gmail.com}{\tt uk.laszlo.toth@gmail.com}
	\end{center}
	
	\vskip .2 in
	
	\begin{abstract}	
		We study the set $\mathcal{S}$ of odd positive integers $n$ with the property ${2n}/{\sigma(n)} - 1 = 1/x$, for positive integer $x$, i.e., the set that relates to odd perfect and odd ''spoof perfect'' numbers. As a consequence, we find that if $D=pq$ denotes a spoof odd perfect number other than Descartes' example, with pseudo-prime factor $p$, then $q>10^{12}$. Furthermore, we find irregularities in the ending digits of integers $n\in\mathcal{S}$ and study aspects of its density, leading us to conjecture that the amount of numbers in $\mathcal{S}$ below $k$ is $\sim10 \log(k)$.
	\end{abstract}

	\section{Introduction}
	
	Let $\mathcal{S}$ denote the set of odd positive integers $n$ with the property
	$$
	\frac{2n}{\sigma(n)} - 1 = \frac1x,
	$$
	where $\sigma$ denotes the sum-of-divisors function and $x$ is a positive integer. These numbers are interesting for several reasons. For instance, if $x$ is prime, $nx$ is an odd perfect number. No such number is currently known, and the abundant amount of restrictions for an odd integer to be perfect (such as those enumerated by Voight \cite{Voight03} and Nielsen \cite{Nielsen07}) suggest that these numbers are either extremely rare or do not exist. 
	
	On the other hand, if $x$ is odd but not prime, the number $nx$ is an \textit{odd spoof perfect number}, i.e., an odd number that would be perfect if only $x$ was prime. These numbers, also referred to as \textit{Descartes numbers} after the discoverer of the only currently known member,
	$$
	D = 198585576189,
	$$
	for which we have $n=9018009$ and $x=22021$, have been subject to considerable research in light of their connection with odd perfect numbers. Despite this, not many things are known about such numbers; for instance, it is not even known if there are an infinite number of positive integers with this property.
	
	A few results however exist. Banks, G\"{u}lo\u{g}lu, Nevans and Saidak \cite{Banks06} showed in 2006 that the only cube-free odd spoof perfect number with fewer than seven distinct prime divisors is Descartes' example. Furthermore, they showed that cube-free spoof odd perfect numbers not divisible by 3 have over a million distinct prime divisors.
		
	Then in 2014, Dittmer \cite{Dittmer13} gave a formal definition of spoof odd perfect numbers and showed that Descartes' example is the only spoof with less than seven distinct quasi-prime factors. In doing so, he defined a \textit{quasi-prime factorization} of a positive integer as a set of pairs $X=\{(x_1,\alpha_1),\ldots,(x_k,\alpha_k)\}$ such that 
	$$
	n=\prod_{i=1}^{k} x_i^{\alpha_i},
	$$
	where $x_i\geq2$ for all $i$, but without the condition that the $x_i$ must be relatively prime. He then used this factorization to define the \textit{spoof $\sigma$-function}, $\tilde{\sigma}$, whose role is analogous to $\sigma$ but is related to $X$ instead of the prime factorization. Under such a factorization, we can say that $n$ is \textit{spoof perfect} if $\tilde{\sigma}(X)=2n$ (although the notation used by Dittmer \cite{Dittmer13} is slightly different). 
	
	We can thus also define ''\textit{spoof abundant}'' and ''\textit{spoof deficient}'' numbers as those having the property $\tilde{\sigma}(X) > 2n$ and $\tilde{\sigma}(X) < 2n$, respectively, under the corresponding quasi-prime factorization $X$.

	\subsection{Scope of this paper}
	
	Instead of considering the number of quasi-prime factors in an odd spoof perfect number $D=pq$, like Dittmer did, we find a new lower bound for $q$ by examining the properties of the set $\mathcal{S}$. In particular, our results show that:
	\begin{theorem}
	Let $D$ denote an odd spoof perfect number such that $D=pq$, where $q\in\mathcal{S}$ and $p$ is the pseudo-prime factor. Furthermore, $D\neq198585576189$. Then $q>10^{12}$.
	\end{theorem}
	In other words, the non-pseudo-prime component of a spoof odd perfect number (other than Descartes' example) is greater than $10^{12}$.
		
	The OEIS sequence \seqnum{A222263} already contains the first 500 terms of $\mathcal{S}$. Using our results presented in this paper, we extend this sequence up to $10^{12}$ (one trillion), then use this dataset to examine various properties of $S$. 
	
	In the remainder of this paper, let $\pi_{\mathcal{S}}(n)$ denote the number of elements in $\mathcal{S}$ up to and including $n$. We define the asymptotic density of the set $\mathcal{S}$ as $\displaystyle \lim \inf \frac{\pi_{\mathcal{S}}(n)}{n}$ and its \textit{Schnirelmann density} as the greatest lower bound of $\displaystyle  \frac{\pi_{\mathcal{S}}(n)}{n}$. We examine these densities in the sections below and submit a few conjectures on their values.

	\section{Computational methods}
	
	The computations were performed on a 6-core Intel i7-7800X processor @ 3.50GHz and took approximately $8$ months to complete. Our algorithm checked each odd $n$ to with the property $\displaystyle \frac{2n}{\sigma(n)} - 1 = \frac1x$ for some positive integer $x$, up to $n=10^{12}$. An outline of the algorithm is shown below.
	
	\begin{algorithm}[h!]
		\caption{Finding positive integers $n\in\mathcal{S}$ smaller than $k$}\label{algo}
		
		\begin{algorithmic}[0]
			
			\State \textbf{Input:} $k$
			\State \textbf{Output:} $n\in\mathcal{S}$ smaller than $k$
			
			\State $results \gets \{\emptyset\}$
			
			\For{\textbf{all odd} $i \leq k$}
			
			\State $m \gets $ \texttt{DivisorSigma}$[1,i] / 2n$
			\State $num \gets $ \texttt{Numerator}$[m]$
			\State $den \gets $ \texttt{Denominator}$[m]$
			\State $diff \gets den - num$
			
			\If{$diff=1$}
				\State $results \gets results \cup \{i\}$
				
				\If{$num\equiv 0 \ (\rm{mod}\ 2)$}
					\State \textbf{print} 'Even spoof perfect number found: $i$'
				\Else{}
					\State \textbf{print} 'Odd spoof perfect number found: $i$'
				\EndIf
			\EndIf			
			\EndFor
			
			\State \Return $results$
			
		\end{algorithmic}
	\end{algorithm}

	This algorithm was executed using Wolfram Mathematica 11.1 and used the built-in \texttt{DivisorSigma[]} function to compute the sum of digits for each candidate.

	Note that this algorithm is naive as it makes no assumptions on the admissibility of a candidate before processing it. However, in the absence of a sufficiently developed theoretical framework there are no practical alternatives to this, to the best of this author's knowledge.
	
	\section{Results}
	
	Our computations revealed no odd spoof perfect number $pq$ other than Descartes' example up to $q=10^{12}$. On the other hand, we found many more even spoof perfect numbers than listed in the OEIS sequence \seqnum{A222263}. Note that these correspond to an even pseudo-prime factor $p$. Table \ref{table-res1} shows the distribution of these numbers within intervals of size $10^k$, for $1 \leq k \leq 12$.
	
	\begin{table}[h!] \caption{Distribution of integers $n\in\mathcal{S}$ in intervals of size $10^k$, $1 \leq k \leq 12$} \label{table-res1}
		\centering
		\begin{tabular}{|l|l|l|}\hline 
			$k$  & $\pi_{\mathcal{S}}(10^k)$ & $\pi_{\mathcal{S}}(10^k)-\pi_{\mathcal{S}}(10^{k-1})$ \\ \hline
			$1$  &  $2$   &  $2$                                          \\ \hline
			$2$  &  $3$   &  $1$                                          \\ \hline
			$3$  &  $7$   &  $4$                                          \\ \hline
			$4$  &  $15$  &  $8$                                          \\ \hline
			$5$  &  $28$  &  $13$                                          \\ \hline
			$6$  &  $48$  &  $20$                                          \\ \hline
			$7$  &  $81$  &  $33$                                          \\ \hline
			$8$  &  $143$ &  $62$                                          \\ \hline
			$9$  &  $227$ &  $84$                                          \\ \hline
			${10}$ &$319$ &  $92$                                          \\ \hline
			${11}$ &$459$ &  $140$                                          \\ \hline
			${12}$ &$692$ &  $233$                                          \\ \hline
		\end{tabular}
	\end{table}
	
	Before examining the density of $\mathcal{S}$ through the lens of our result set, we first take a look at some of its characteristics. We begin with the congruence classes formed by the $n\in\mathcal{S}$, as shown in Table \ref{table-res2}. Our results seem to indicate that the $n$ are uniformly distributed into residue classes mod $8$.
	\begin{table}[h!] \caption{Congruence of integers $n\in\mathcal{S}$ up to $10^{12}$} \label{table-res2}
		\centering
		\begin{tabular}{|l|l|}\hline 
			Residue class  & Amount of integers $n\in\mathcal{S}$ in the different residue classes \\ \hline
			$1$ mod $8$  &  $149$                                             \\ \hline
			$3$ mod $8$  &  $183$                                             \\ \hline
			$5$ mod $8$  &  $175$                                             \\ \hline
			$7$ mod $8$  &  $185$                                            \\ \hline
		\end{tabular}
	\end{table}

	In hopes of finding an irregularity among the members of our data set, we examined the distribution of their ending digits. And indeed, our results show a strong bias in favour of numbers ending in $5$ than those ending in other digits. The distribution of ending digits is shown in Table \ref{table-res3}.
	\begin{table}[h!] \caption{Distribution of ending digits of integers $n\in\mathcal{S}$ up to $10^{12}$} \label{table-res3}
		\centering
		\begin{tabular}{|l|l|}\hline 
			$l$  & Amount of $n\in\mathcal{S}$ with $l$ as last digit \\ \hline
			$1$  &  $51$                                             \\ \hline
			$3$  &  $46$                                             \\ \hline
			$5$  &  $492$                                             \\ \hline
			$7$  &  $54$                                            \\ \hline
			$9$  &  $49$                                            \\ \hline
		\end{tabular}
	\end{table}

	\subsection{Density}
		
	We pursue our analysis by studying the density of $\mathcal{S}$, i.e., the ratio $\displaystyle \frac{\pi_{\mathcal{S}}(n)}{n}$ for positive integer $n$. We have plotted this density at each $n\in\mathcal{S}$ on a log-log graph, shown in Figure \ref{fig-density} below.
	\begin{figure}[h!] 
		\caption{Density of $\mathcal{S}$ up to $k=10^{12}$}
		\label{fig-density}
		\centering
		\includegraphics[]{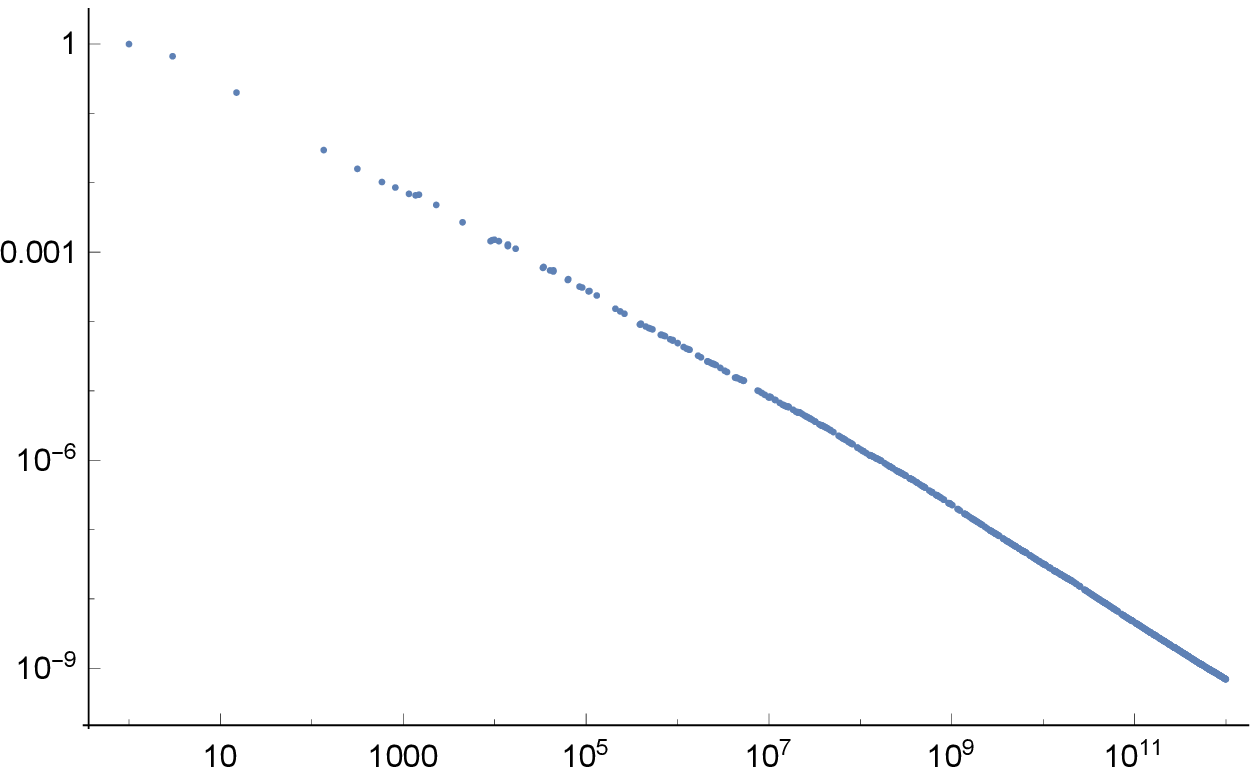}
	\end{figure}
	
\begin{figure}[ht!]
		\caption{Density of $\mathcal{S}$ (blue) and $A(k)=10 \log(k)/k$ (orange) up to $k=10^{12}$}
		\label{fig-density2}
		\centering
		\begin{minipage}[b]{0.45\linewidth}
			\includegraphics[scale=0.65]{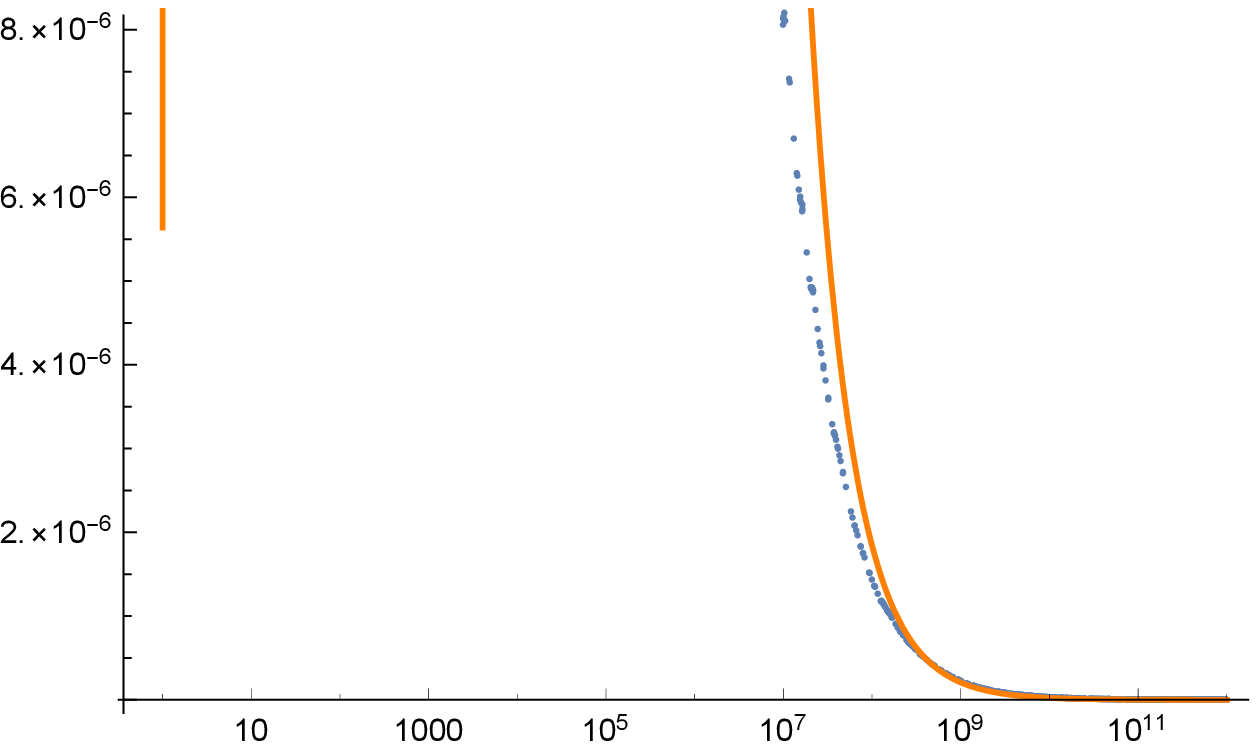}
			\label{fig:minipage1}
		\end{minipage}
		\quad
		\begin{minipage}[b]{0.45\linewidth}
			\includegraphics[scale=0.65]{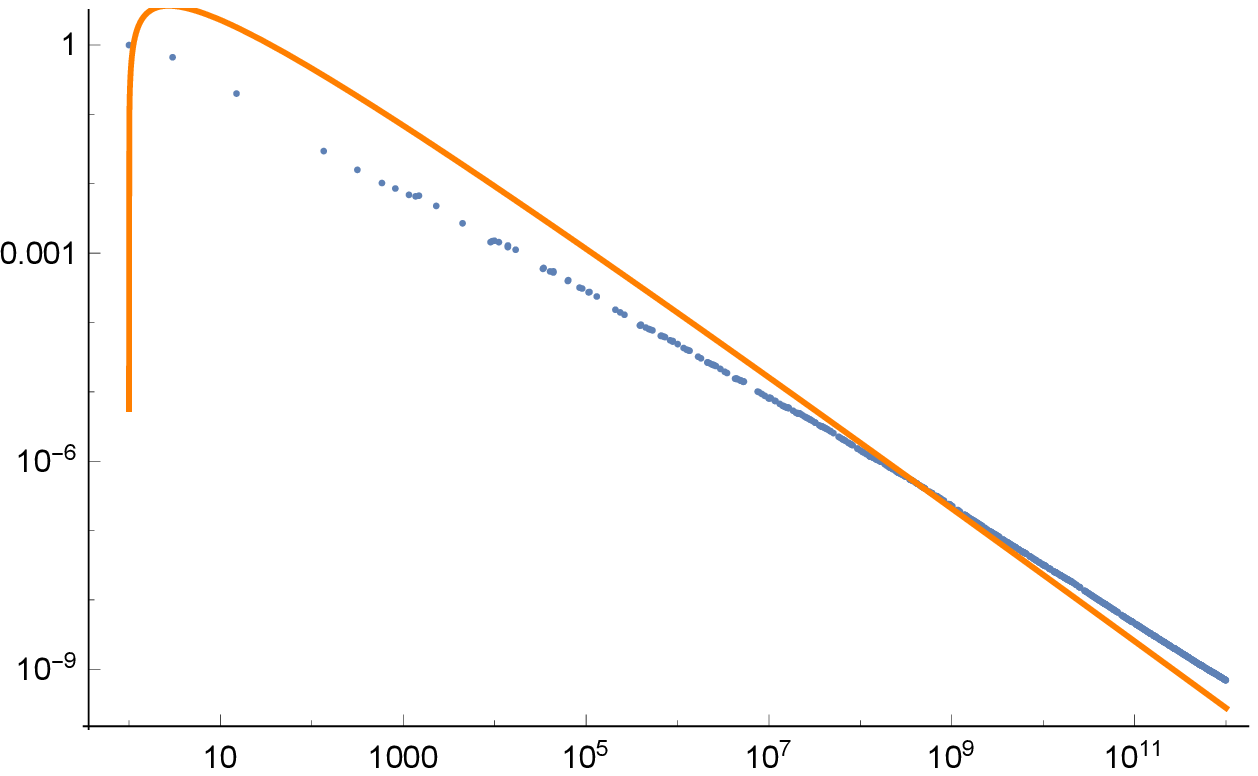}
			\label{fig:minipage2}
		\end{minipage}
	\end{figure}

	This suggests that the density follows a curve of the type 
	$$
	A(k) = \frac{\alpha \log(k)}{k},
	$$
	with real $\alpha > 0$. We have found that a value around $\alpha=10$ provides a good fit to our experimental data, which we show in Figure \ref{fig-density2} below on a log-linear and log-log graph.
	
	If this is true, it would naturally follow that:
	\begin{conjecture}
		$\pi_{\mathcal{S}}(n) \sim 10 \log(n)$.
	\end{conjecture}	
	These results suggest several interesting properties of $\mathcal{S}$. On the one hand, we are tempted to conjecture that:
	\begin{conjecture} \label{conj-density-zero}
		The asymptotic density of $\mathcal{S}$ is $0$.
	\end{conjecture}
	On the other hand, even though $\{1\}\in\mathcal{S}$, our data suggests that the fraction $\displaystyle \frac{\pi_{\mathcal{S}}(k)}{k}$ tends to $0$ as $k\to\infty$, and so we also conjecture that:
	\begin{conjecture}
		The Schnirelmann density of the set $\mathcal{S}$ is $0$.
	\end{conjecture}

	\section{Conclusion and further work}

	Nothing in our results suggests that another spoof perfect number exists, or even that there are an infinite number of positive integers in $\mathcal{S}$. On the other hand, we made a not entirely unreasonable conjecture that the density of $\mathcal{S}$ is $0$, providing further evidence of the scarcity of such numbers.
	
	There are several ways to extend the results in this paper. First, it is easy to continue computations above $10^{12}$ with sufficient computing resources. Furthermore, one might be inspired to write a more clever algorithm by taking a look at the candidate integers before computing their sum-of-divisors function. This operation is expensive and could be avoided if we have additional information about such numbers. However, this is not currently the case, making it difficult to escape combing through odd integers up to a given number. Perhaps the considerable effort already expended on similar questions regarding odd perfect numbers could be adapted to spoof odd perfect numbers as well.

	Classification: 11A51, 11N25, 11B83

\end{document}